\documentclass{article}
\usepackage[english]{babel}
\usepackage{amsmath,amsfonts,amsthm,amssymb,amscd}

 \usepackage{verbatim}

\renewcommand{\Re}{\mathop{\rm Re}\nolimits}
\renewcommand{\Im}{\mathop{\rm Im}\nolimits}

\newcommand{\p}{\partial}
\newcommand{\e}{\varepsilon}

\newcommand{\C}{{\mathbb C}}
\newcommand{\R}{{\mathbb R}}

\newcommand{\pP}{{\mathbb P}}
\newcommand{\I}{{\mathbb I}}
\newcommand{\E}{{\mathbb E}}

\newcommand{\N}{{\mathbb N}}
\newcommand{\la}{\lambda}

\newcommand{\ty}{\infty}

\newcommand{\aA}{{\cal A}}

\newcommand{\DD}{{\cal D}}

\newcommand{\FF}{{\cal F}}
\newcommand{\GG}{{\cal G}}

\newcommand{\KK}{{\cal K}}

\newcommand{\OO}{{\cal O}}

\newcommand{\UU}{{\cal U}}
\newcommand{\VV}{{\cal V}}

\newcommand{\lag}{\langle}
\newcommand{\rag}{\rangle}

\newcommand{\dd}{{\textup d}}

\newcommand{\supp}{\mathop{\rm supp}\nolimits}

\theoremstyle{plain}
\newtheorem{theorem}{Theorem}[section]
\newtheorem{lemma}[theorem]{Lemma}
\newtheorem{proposition}[theorem]{Proposition}

\newtheorem{condition}[theorem]{Condition}
\theoremstyle{remark}
\newtheorem{remark}[theorem]{Remark}
\newtheorem{example}[theorem]{Example}
\newcommand{\de}{\delta}
\numberwithin{equation}{section}
\begin{document}
\author{Vahagn Nersesyan}
\date{}
\title{ Growth of Sobolev norms and controllability of  Schr\"odinger equation}

 \maketitle

\begin{center}
 Laboratoire de Math\'ematiques,
Universit\'e de Paris-Sud XI\\ B\^atiment 425, 91405 Orsay Cedex,
France\\ E-mail: Vahagn.Nersesyan@math.u-psud.fr
\end{center}
{\small\textbf{Abstract.} In this paper  we obtain a stabilization
result for the   Schr\"odin\-ger equation under
generic assumptions on the potential. 
Then we consider the  Schr\"odin\-ger equation with a potential
which has a random time-dependent amplitude. We show that if the
distribution of the amplitude  is sufficiently non-degenerate,
then any trajectory of system is almost surely non-bounded in
Sobolev spaces.
  }\\\\
 \tableofcontents

\section{Introduction}\label{S:intr}
We consider the problem
\begin{align}
i\dot z   &= -\Delta z+V(x)z+ u(t) Q(x)z,\,\,\,\,x\in D,\label{E:hav1}\\
z\arrowvert_{\partial D}&=0,\label{E:ep1}\\
 z(0,x)&=z_0(x),\label{E:sp1}
\end{align} where $D\subset\R^m$ is a bounded domain with smooth boundary, $V,Q\in C^\ty(\overline{D},\R)$ are  given functions,  $u$ is the control, and $z$
is the state. Under some hypotheses on $V$ and $Q$ (see Condition
\ref{C:p}), we prove a stabilization result for   problem
(\ref{E:hav1}), (\ref{E:ep1}). Then this result is applied to show
that almost any trajectory of random  Schr\"odinger equation is
non-bounded in Sobolev spaces.   As is shown in Section
\ref{S:gen}, the hypotheses on $V$ and $Q$ are in a sense generic.

Let us recall some previous results on  controllability of
Schr\"odinger equation. A general negative result for bilinear
control systems is obtained by  Ball,  Marsden and Slemrod
\cite{BMS}. Application of this result to (\ref{E:hav1}),
(\ref{E:ep1}) implies that the set of attainable points from any
initial data in $H^2$ admits a dense complement in $H^2$. We refer
the reader to the papers \cite{RSDR,TR,ALT,ALAL,AC} and the
references therein for controllability of finite-dimensional
systems. In \cite{BCH}, Beauchard  proves that exact
controllability result is possible to obtain if one chooses
properly the phase space. More precisely, in the case $m=1,V(x)=0$
and $Q(x)=x$    exact controllability of the problem is proved  in
$H^7$-neighborhood of the  eigenstates. A stabilization property
for finite-dimensional approximations of Schr\"odinger equation is
obtained by Beauchard et al., in \cite{BCMR}, which was later
generalized by Beauchard and Mirrahimi \cite{BM} to the
infinite-dimensional  case for $m=1,V(x)=0$ and $Q(x)=x$ (see also
the paper by Mirrahimi \cite{Mir}). Recently Chambrion et al.
\cite{CH}, under some assumptions on $V,Q\in
C^\ty(\overline{D},\R)$, derived the approximate controllability
of    (\ref{E:hav1}), (\ref{E:ep1}) in $L^2$ from the
controllability of finite-dimensional projections. See also the
papers \cite{L,BQ,MZ,BP,Z,DGL} and the references therein for
controllability results by boundary controls and controls
supported in a given subdomain and  the book  \cite{C} by Coron
for introduction to the later developments and methods in the
control theory  of nonlinear systems.

The main result of this paper states that any neighborhood of the
first eigenfunction of operator $-\Delta +V$ is attainable from
any initial point $z_0\in H^2 $. This result, combined with the
time reversibility property of the system and the fact that the
equation is linear, implies  approximate controllability
property in~$L^2$.

Let us describe in a few words the main ideas of the proof. As $V,
Q$ and~$u$ are real-valued, the $L^2$ norm is preserved by the
flow of the system. Thus it suffices to consider the restriction
of  (\ref{E:hav1}), (\ref{E:ep1}) to the unit sphere $S$ in $L^2$.
We introduce a Lyapunov  function $\VV(z)$ that controls the
$H^2$-norm of $z$. The infimum of $\VV$ on the sphere $S$ is
attained at the first eigenfunction $e_{1,V}$ of the operator
$-\Delta+V$. Using the idea of generating trajectories with
Lyapunov techniques from \cite{BCMR}, we choose a feedback law
$u(z)$ such that the function $\VV$ decreases on the solutions of
the corresponding system:
$$
\VV(\UU_{t}(z_0,u))<\VV(z_0),\,\,\,t>0,
$$where $\UU_{t }(\cdot,u)$ is the resolving operator of  (\ref{E:hav1}), (\ref{E:ep1}).
Then iterating this construction and using the fact that the
system is autonomous, we prove that the $H^2$-weak $\omega$-limit
set of any solution contains   the minimum point of function
$\VV$, i.e. the eigenfunction $e_{1,V}$ (see Sections \ref{S:STAB}
and \ref{S:??}).

The ideas of the proof  work  also in the case of nonlinear
equation. We think that the result holds also in the spaces $H^l,
l>2$. This case will be treated in a later paper.

We next use  the above-mentioned controllability result  to study
the large time behavior of  solutions of  random Schr\"odinger
equation. We show that  if the distribution of the random
potential is sufficiently non-degenerate (see Condition
\ref{C:c}), then the trajectories of the system are almost surely
non-bounded. It is interesting to compare this result with that of
Eliasson and Kuksin \cite{EK}, where KAM-technique is applied to
prove the reducibility of a linear Schr\"o\-din\-ger equation with
time-quasiperiodic potential. In particular, it is proved that for
most   values of the frequency vector the Sobolev norms of the
solutions are bounded. Examples of non-bounded solutions of 1D
linear Schr\"odinger equation with some random potentials are
constructed in \cite{BG,EKS}, where also the growth rate estimates
are given. Our assumptions  on the distribution of the potential
are more general, and the proof works also in the case of
nonlinear equation. However, at this level of generality, we do
not have any lower bound on the rate of growth of Sobolev norms.

The idea of the proof is to show that the first entrance time to
any ball  centered at the origin in $H^{-\e}$ is almost surely
finite. This implies immediately that  almost any trajectory of
the system approaches the origin arbitrarily closely in $H^{-\e}$.
Combining this with the fact that the $L^2$-norm is preserved, we
conclude that almost any trajectory is non-bounded in $H^l$ for
any $l>0$.

In conclusion, let us note that the results of this paper imply
the irreducibility in $L^2$ of the Markov chain associated with
(\ref{E:hav1}).   This property is not sufficient to prove the
ergodicity of the dynamics generated by the Schr\"odinger equation
with random potential. However, in the case of finite-dimensional
approximations, that question is treated in the paper \cite{N}, in
which an exponential mixing property is established. We hope the
methods developed in this work will help to tackle the
infinite-dimensional case.

\textbf{Acknowledgments.}  The author would like to thank  Armen
Shirikyan for his guidance and encouragements.

\vspace{39pt} \textbf{Notation}\\\\ In this paper we use the
following notation. Let $D \subset \R^m, m\ge1$ be a bounded
domain with smooth boundary. Let  $H^s:=H^s(D)$  be the Sobolev
space of order $s\in\R$ endowed with the norm $\|\cdot\|_s$.
Consider the operators $ -\Delta z +Vz,z\in\DD(- \Delta +V):=H_0^1
\cap H^2, $ where $V\in C^\ty(\overline{D},\R)$. We denote by
$\{\la_{j,V} \}$ and $\{e_{j,V} \}$ the sets of eigenvalues and
normalized eigenfunctions  of $
-\Delta+V$.  
Let $\langle\cdot,\cdot\rangle$ and $\|\cdot\|$ be the scalar
product and the norm in the space $L^2 $. Let $S$ be the unit
sphere in $L^2 $. For a Banach space $X$,  we shall denote by
$B_X(a,r)$ the open ball of radius $r>0$ centered at $a\in X$. The
set of irrational numbers is denoted by~$\I$.

\section{Preliminaries}\label{S:PREM}

The following lemma shows the well-posedness of system
(\ref{E:hav1})-(\ref{E:sp1}).
\begin{lemma}\label{L:LD}

For any   $z_0\in H_0^1\cap H^2$ $(z_0\in L^2)$ and for any
  $u\in L^1_{loc}([0,\ty),\R)$
  problem
(\ref{E:hav1})-(\ref{E:sp1}) has a unique solution $z\in
C([0,\infty), H^2)$ $(z\in C([0,\infty),L^2))$. 
Furthermore, the resolving operator $\UU_t(\cdot,u):L^2
\rightarrow L^2 $ taking $z_0$ to
    $z(t)$ satisfies the relation
\begin{align}
\|\UU_t(z_0,u)\|&=\|z_0\|,\,\,\,t\ge0.\label{E:barev}
\end{align}
\end{lemma}
See  \cite{CW} for the proof. 
 Notice that the
conservation of $L^2$-norm   implies that it suffices to consider
the  controllability properties of (\ref{E:hav1}), (\ref{E:ep1})
on the unit sphere $S$.

In Section \ref{intro}, we replace the control $u$ by a random
process. Namely,  we consider the equation

\begin{equation*}
i\dot z = -  \Delta z+V(x)z +\beta(t) Q(x)z,\,\,\,\,x\in D,
\end{equation*}where   $\beta(t)$ is a random process of the
form
\begin{equation}\label{E:pu}
\beta (t)=\sum_{k=0}^{+\infty}I_{k}(t)\eta_k(t-k),\,\,\,\,\,t\ge
0.
\end{equation}
Here $I_{k}(\cdot)$ is the indicator function of the interval
$[k,k+1)$ and $\eta_k$ are independent identically distributed
(i.i.d.) random variables in $L^2([0,1],\R)$.

Let $z_0$ be    $L^2$-valued  random variable independent of
$\{\eta_k\}$. Denote by $\FF_k$  the $\sigma$-algebra generated by
$z_0,\eta_0,\dots,\eta_{k-1}$.

\begin{lemma}\label{L:MP}
Under above conditions,  $\UU_k(\cdot,\beta)$ is a homogeneous
Markov chain with respect to $\FF_k$.
\end{lemma}
This lemma is proved by standard arguments (e.g., see \cite{OKS}).

\section{Controllability of the Schr\"odinger equation}\label{S:er}

\subsection{Stabilization result }\label{S:STAB}

Let us  introduce the Lyapunov function
\begin{equation*}
\VV(z): =\alpha \|( -\Delta+{V} )P_{1,V}z\| ^2+1-|\lag z,e_{1,V}
\rag |^2,\,\,\,\, z\in S\cap H^1_{0}\cap H^2,
\end{equation*}where
 $\alpha>0$
  and   $P_{1,V}z:=z-\lag z,e_{1,V}\rag e_{1,V}$ is  the
orthogonal projection in $L^2$ onto the closure of the vector span
of $\{e_{k,V}\}_{k\ge 2}$. Notice that $\VV(z)\ge0$ for all $ z\in
S\cap H^1_{0}\cap H^2$ and $\VV(z)=0$ if and only if $z=ce_{1,V},
|c|=1$. For any $ z\in S\cap H^1_{0}\cap H^2$, we have
\begin{align*}
\VV(z)\ge\alpha \|( -\Delta+{V} )P_{1,V}z\| ^2\ge\frac{\alpha}{2}
\|  \Delta (P_{1,V}z)\| ^2-C_1\ge\frac{\alpha}{4} \|  \Delta z\|
^2-C_2.
\end{align*}Thus
\begin{equation}\label{E:sah}
C(1+\VV(z))\ge \|z\|_2
\end{equation}for some constant $C>0$.
Following the ideas of \cite{BCMR}, we wish to choose a feedback
law $u (\cdot)$ such that
\begin{equation*}
\frac{\dd}{\dd t}\VV(z(t))\le0
\end{equation*}for the  solution
$z(t)$ of (\ref{E:hav1})-(\ref{E:sp1}). Let us assume that $\Delta
z(t)\in H_0^1\cap H^2$ for all $t\ge0.$
  Using (\ref{E:hav1}), we get
 \begin{align*}
\frac{\dd}{\dd t}\VV(z(t))&=2\alpha\Re(\langle (-\Delta +V)
P_{1,V}{\dot{
z}}, (-\Delta +V)  P_{1,V}z\rangle)-2\Re(\lag\dot z,e_{1,V}\rag\lag e_{1,V},z\rag)\nonumber\\
&= 2\alpha\Re(\langle  (-\Delta +V) P_{1,V}(i\Delta z -iVz-iu  Qz
), (-\Delta
+V)  P_{1,V}z\rangle)\nonumber\\
&\quad-2\Re(\lag  i\Delta z -iVz-iuQz,e_{1,V}\rag\lag
e_{1,V},z\rag).
\end{align*}Integrating by parts and using the fact that $$(-\Delta +V)  P_{1,V}z|_{\p D}=z|_{\p D}=e_{1,V}|_{\p D}=0,$$ we obtain
 \begin{align*}
 & 2\alpha\Re(\langle  -i(-\Delta +V)^2  P_{1,V}z, (-\Delta +V)  P_{1,V}z\rangle) -2\Re(\lag  i\Delta z -iVz ,e_{1,V}\rag\lag e_{1,V},z\rag)\nonumber\\
&=  2\alpha\Re(\langle   i\nabla (-\Delta +V)   P_{1,V}z, \nabla
(-\Delta +V) P_{1,V}z\rangle)\nonumber\\&\quad+2\alpha\Re(\langle
-iV(-\Delta +V) P_{1,V}z, (-\Delta +V)
P_{1,V}z\rangle)\nonumber\\&\quad+2\la_{1,V}\Re(\lag  i z
,e_{1,V}\rag\lag e_{1,V},z\rag)=0.
\end{align*}
Thus
\begin{align*}
\frac{\dd}{\dd t}\VV(z(t)) = 2u\Im(\alpha\lag (-\Delta +V)
P_{1,V}(Qz), (-\Delta +V) P_{1,V}z\rangle -\lag Qz,e_{1,V}\rag\lag
e_{1,V},z\rag).
\end{align*}
Let us  take
\begin{equation}\label{E:fed3}
  {{ {u}}}(z):=-\delta\Im(\langle \alpha (-\Delta +V)   P_{1,V}(Qz ), (-\Delta
 +V)
 P_{1,V}z\rangle
 - \lag Qz,e_{1,V}\rag\lag e_{1,V},z\rag),
 \end{equation} where $\de>0$ is a small constant. Then
\begin{equation}\label{E:fed35}
 \frac{\dd}{\dd t} \VV(z(t))=-\frac{2}{\delta}u^2(z(t)).
\end{equation}Consider the equation
\begin{equation}
i\dot z    = -\Delta z+V(x)z+ u(z) Q(x)z,\,\,\,\,x\in
D.\label{E:hav22}
\end{equation}
 \begin{proposition}\label{T:abc}For any $z_0\in H_0^1\cap H^2$   problem
 (\ref{E:hav22}), (\ref{E:ep1}), (\ref{E:sp1}) has a unique
 solution $z\in C([0,\ty),  H_0^1\cap H^2)$.  Moreover, the
   following properties hold.
\begin{enumerate}
\item [(i)] If $\Delta z_0\in H_0^1\cap H^2$, then $\Delta z\in C([0,\ty),  H_0^1\cap
H^2)$.
\item  [(ii)] Let $\UU_t(\cdot): H_0^1\cap H^2\rightarrow
   H_0^1\cap H^2$ be the resolving operator. If  $T>0$, $z_n\in  H_0^1\cap H^2$ and $z_n\rightharpoonup
 z_0$ in $H^2$, then $\UU_{T}(z_{n_k})\rightharpoonup\UU_T(z_0)$ in $H^2$ for      some sequence $k_n\ge1$.
\end{enumerate}
  \end{proposition}
  \begin{proof}[Sketch of the proof] The local well-posedness of (\ref{E:hav22}), (\ref{E:ep1}) and
  (\ref{E:sp1}) is standard (see \cite{CW}). From the
  construction of the feedback $u$ it follows that a finite-time
  blow-up is impossible. Hence the solution is global in time. To
  prove the rest of the theorem, it suffices to show that $u(z_n)\rightarrow
  u(z_0)$ for any $z_n\in  H_0^1\cap H^2$ such that $z_n\rightharpoonup
 z_0$ in $H^2$. Notice that (\ref{E:fed3}) and the fact that $Q$ is real imply  that
 \begin{align*}
u(z)=-\Im(\langle \alpha Q(-\Delta +V)     z  , (-\Delta
 +V)
  z\rangle +\tilde{u}(z)=\tilde{u}(z),
 \end{align*}where $\tilde{u}(z_n)\rightarrow
  \tilde{u}(z_0).$ This completes the proof.
  \end{proof}
Thus if $z_0,\Delta z_0\in H_0^1\cap H^2$, then (\ref{E:fed35}) is
verified for $z(t)=\UU_t(z_0)$. A density argument proves the
identity for any $z_0\in H_0^1\cap H^2$.

  Let us assume that the functions $V$ and $Q$ satisfy the following
condition.

\begin{condition}\label{C:p}
The functions $V,Q\in C^\ty(\overline{D},\R) $ are such that:
\begin{enumerate}
\item [(i)] $\langle Qe_{1,V},e_{j,V}\rangle\neq0$
  for all $j\ge1$,
\item  [(ii)]$\la_{1,V}-\la_{j,V}\neq \la_{p,V}-\la_{q,V}$ for all  $j,p,q\ge1$ such that
  $\{1,j\}\neq\{p,q\}$ and $j\neq 1$.

\end{enumerate}
\end{condition}
The below theorem is the main result of this section.
 \begin{theorem}\label{T:stab}
 Under Condition \ref{C:p},  there is a finite  or countable set $J\subset \R_{+}^*$ such that for any
  $\alpha\notin J$ and $z_0\in S\cap H^1_{0}\cap H^2$ with $\lag z_0,e_{1,V} \rag\neq0$
  and $0<\VV(z_0)<1$  there is a sequence
 $k_n\ge1$ verifying
 $$
 \UU_{k_n}(z_0 )\rightharpoonup  c e_{1,V}\,\,\,\,\text{in
$H^2$},
 $$where $c\in\C,|c|=1.$
 \end{theorem}See Subsection \ref{S:??} for the proof of this
 theorem. The following lemma shows that the hypothesis on the
 initial condition $z_0$ is not restrictive.
 \begin{lemma}\label{L:eee}
 For any $z_0\in S$ there is a control $u\in C^\ty([0,\ty),\R)$
 and a time $k\ge1$ such that  $\lag\UU_k(z_0,u),e_{1,V}\rag\neq0.$
 \end{lemma}
 \begin{proof}
       It suffices to find a control $u$ and a time $k\ge1$ such that
\begin{equation}\label{E:ite}
\| \UU_k (z_0,u)-ce_{1,V} \|<\sqrt{2}
\end{equation}
for some  $c\in\C, |c|=1$. Take any $\hat {z}_0\in S\cap H^1_0\cap
H^2 $ such that  $\lag \hat z_0,e_{1,V} \rag\neq0$ and
$$\|z_0-\hat z_0\|<\frac{\sqrt 2}{2}.$$ By Theorem \ref{T:stab}, there is a control $u\in C^\ty([0,\ty),\R)$
 and a time $k\ge1$ such that
$$
\| \UU_k (\hat z_0,u)-ce_{1,V} \|<\frac{\sqrt 2}{2}.
$$Using the fact that  the $L^2$-distance between two solutions of (\ref{E:hav1}),
(\ref{E:ep1}) with the same control is constant, we obtain
(\ref{E:ite}).

 \end{proof}
\subsection{Approximate controllability}\label{S:APP}
Before proving      Theorem \ref{T:stab}, let us give an
application of the result.
 For any
$d>0$ define the set
$$
C_{d}=\{u\in C^\ty([0,\ty),\R): \sup_{t\in[0,\ty)} | u(t) |<d \}.
$$
 We say that problem (\ref{E:hav1}), (\ref{E:ep1}) is
approximately controllable in $L^2$  at integer times if for any
$\e,d>0$ and for any points $z_0,z_1\in S$ there is a time $k\in
\N$ and a control $u\in C_{d}$ such that
\begin{equation}
\|\UU_k(z_0,u)-z_1\|<\e.\nonumber
\end{equation}
\begin{theorem}\label{T:con}
Under Condition \ref{C:p},   problem (\ref{E:hav1}), (\ref{E:ep1})
is approximately controllable in $L^2$ at integer times.
\end{theorem}
\begin{proof}

  Theorem  \ref{T:stab} implies that for any $z\in S\cap H^1_{0}\cap H^2 $ there is $u\in C_d$ such that
 \begin{equation} \label{E:opo1} \|\UU_{k}(z,u)-
e_{1,V}\|<\frac{\e}{2}
 \end{equation}for some $k\ge1$. As the
$L^2$-distance between two solutions of (\ref{E:hav1}),
(\ref{E:ep1}) with the same control is constant, by a density
argument,  we get that for any $z\in S$ a control $u\in C_d$
exists such that  (\ref{E:opo1}) holds.

Here we need the following result    often referred as time
reversibility property of  Schr\"odinger equation.
\begin{lemma}\label{L:iop}
Suppose that $\UU_k(\bar{z},w)=\bar{y}$ for some $z\in L^2$, $w\in
C_d$ and $k\ge1$. Then $\UU_k(y,u)=z$, where $u(t)=w(k-t)$.
\end{lemma}
The proof of this lemma is clear.  Let us fix any $z_0,z_1\in S$
and let $u_0,w\in C_{d}$ be such that
\begin{align}
&\|\UU_{k_1}(\bar z_1,w)-e_{1,V}\|<\frac{\e}{2},\nonumber\\
&\|\UU_{k_0}(z_0,u_0)-e_{1,V}\|<\frac{\e}{2}\nonumber
\end{align} for some $k_0,k_1\ge1$. Define  $y:=\overline{\UU_{k_1}(\bar z_1,w)}$. Then
by Lemma \ref{L:iop}, we have $\UU_{k_1}(y,u_1)=z_1$, where
$u_1(t):=w(k_1-t)$. Again using the fact that $L^2$-distance
between two solutions of (\ref{E:hav1}), (\ref{E:ep1}) with the
same control is constant, we get
\begin{equation}
\|\UU_{k_1}(e_{1,V},u_1)-z_1\|=\|e_{1,V}-y\|<\frac{\e}{2}.\nonumber
\end{equation}

 Taking  $k=k_0+k_1$ and $\hat{u}(t)=u_0(t)$, $t\in[0,k_0)$ and
$\hat{u}(t)=u_1(t-k_0)$, $t\in[k_0,\ty)$, we obtain
$$
\|\UU_k(z_0,\hat{u})-z_1\|<\e.
$$
Finally, using the continuity of $\UU_k(z_0,\cdot)$, we find $u\in
C_{d}$ satisfying
$$
\|\UU_k(z_0,{u})-z_1\|<\e.
 $$
\end{proof}
\begin{remark}
We note that for $m=1, Q(x)=x$ a stronger result is obtained by K.
Beauchard and M. Mirrahimi \cite{BM} in the case of the space
$L^2$. They show an approximate stabilization result of
eigenstates. The proof of this result remains literally the same
for system (\ref{E:hav1}), (\ref{E:ep1}) under Condition
\ref{C:p}. One should just  pay attention to the fact that in the
case of any space dimension $m$  the spectral gap property for the
eigenvalues used in \cite{BM} does not hold. The argument can be
replaced by 
Lemma \ref{L:ner}.

\end{remark}

\subsection{Proof of Theorem \ref{T:stab}}\label{S:??}

\vspace{6pt}\textbf{Step 1.} Let us suppose that $u(\UU_t(z_0))=
0$ for all $t\ge0$. Then
\begin{equation}\label{E:opm}
\UU_t(z_0)=\sum_{j=1}^\ty e^{-i\la_{j,V}t}\langle
z_0,e_{j,V}\rangle e_{j,V}.
\end{equation}
   Substituting (\ref{E:opm}) into (\ref{E:fed3}), we get
\begin{align*}
 0&= \sum_{j=1,k=2}^\ty \alpha \la_{k,V}  \langle z_0,e_{j,V} \rangle
\langle e_{k,V} ,z_0\rangle \langle (-\Delta +V)
(P_{1,V}(Qe_{j,V})) ,e_{k,V} \rangle e^{-i(\la_{j,V} -\la_{k,V}
)t} \nonumber\\&\quad- \sum_{j=1,k=2}^\ty \alpha \la_{k,V} \langle
e_{j,V} ,z_0\rangle\langle z_0,e_{k,V} \rangle \langle e_{k,V},
(-\Delta +V) (P_{1,V}(Qe_{j,V})) \rangle e^{i(\la_{j,V}
-\la_{k,V} )t}\nonumber\\
&\quad- \sum_{j=1}^\ty \langle z_0,e_{j,V} \rangle \langle e_{1,V}
,z_0\rangle \langle Qe_{j,V} ,e_{1,V} \rangle e^{i(\la_{1,V}
-\la_{j,V} )t} \nonumber\\&\quad+ \sum_{j=1}^\ty \langle
e_{j,V},z_0 \rangle \langle z_0,e_{1,V} \rangle \langle Qe_{j,V}
,e_{1,V} \rangle e^{-i(\la_{1,V} -\la_{j,V} )t} \end{align*}
\begin{align}\label{E:mmpp}
 &=\sum_{j=2,k=2}^\ty
P(z_0,Q,j,k) e^{-i(\la_{j,V} -\la_{k,V} )t}
\nonumber\\&\quad+\sum_{j=2}^\ty\Big[ (\alpha\la_{j,V}
\langle(-\Delta +V) (P_{1,V}(Qe_{1,V})),e_{j,V}\rangle+\langle
Qe_{j,V} ,e_{1,V}\rangle )\nonumber\\&\quad\quad\quad\times
\langle z_0,e_{1,V} \rangle \langle e_{j,V} ,z_0\rangle
e^{i(\la_{1,V} -\la_{j,V} )t}\Big]
\nonumber\\&\quad-\sum_{j=2}^\ty \Big[(\alpha\la_{j,V}
\langle(-\Delta +V) (P_{1,V}(Qe_{1,V})),e_{j,V}\rangle+ \langle
Qe_{j,V} ,e_{1,V}\rangle ) \nonumber\\&\quad\quad\quad
\times\langle e_{1,V},z_0 \rangle \langle z_0,e_{j,V} \rangle
e^{-i(\la_{1,V} -\la_{j,V} )t}\Big],
\end{align}where $P(z_0,Q,j,k)$ is a constant.
 In view of Condition \ref{C:p}, (ii),
  Lemma \ref{L:ner} below implies that the coefficients of
exponential functions in (\ref{E:mmpp}) vanish. Condition~\ref{C:p},~(i),  implies that the set 
\begin{align*}J:=\{\alpha\in\R:\alpha\la_{j,V}  \langle(-\Delta +V)
(P_{1,V}(Qe_{j,V})),e_{1,V}\rangle+ \langle Qe_{j,V}&
,e_{1,V}\rangle=0 \nonumber\\&\text{for some
$j\ge1$}\}\end{align*} is       finite  or countable. Thus we get
that $z_0=c e_{1,V}$ for some $c\in\C, |c|=1$ which is a
contradiction to $\VV(z_0)>0$.
Thus  there is a time $t_0>0$ such that $ u(\UU_{t_0}(z_0))\neq0$
and
$$
\VV(\UU_k(z_0))-\VV(z_0)=-\frac{2}{\de}\int_0^k u^2(\UU_s(z_0))\dd
s<0
$$for any $k\ge t_0$.

\vspace{6pt}\textbf{Step 2.} Let $\KK$ be the $H^2$-weak
$\omega$-limit set of the trajectory for (\ref{E:hav22}),
(\ref{E:ep1}) issued from $z_0$, i.e.
\begin{align*}
\KK=\{z\in   H_{0}^{1}\cap H^2& :  \UU_{k_n}(z_0)
 \rightharpoonup z\,\,\,\text{in $H^{2  }$}
    \text{ for some  }   k_n\rightarrow\ty \}.
\end{align*}Let
$$
m:=\inf_{z\in\KK}\VV(z).
$$
This infimum is attained, i.e. there is $e\in \KK$ such that
\begin{equation*}
\VV(e)=\inf_{z\in\KK}\VV(z).
\end{equation*}Indeed, take any minimizing sequence $z_n\in\KK$, so that $\VV(z_n)\rightarrow
m$. By (\ref{E:sah}), $z_n$ is bounded in $H^{2}$. Thus, without
loss of generality, we can assume that $z_n\rightharpoonup e$ in
$H^2$. This implies that
$\VV(e)\le\liminf_{n\rightarrow\ty}\VV(z_n)= m$.
 Let us show that
$e\in \KK$. We can choose  a sequence $k_n\ge1$ such that
\begin{equation}\label{E:ppp}
\|\UU_{k_n}(z_0 )-z_n\| \le\frac{1}{n}.
\end{equation}As $\UU_{k_n}(z_0)$ is bounded in $H^2$, without
loss of generality, we can suppose that  $\UU_{k_n}(z_0
)\rightharpoonup \tilde{e}$, $\tilde{e}\in S\cap H_0^1\cap H^2$.
Clearly, (\ref{E:ppp}) implies that $e=\tilde{e}$, hence $e\in\KK$
and  $\VV(e)= m$.

Let us
 show that $\VV(e)=0$. Suppose that $\VV(e)>0$. As $\VV(e)\le\VV(z_0)<1,$ we have $\lag e,e_{1,V}\rag\neq0.$ Then, by Step 1, there is a  time $k\ge1$
 such that $\VV(\UU_k (e))<\VV(e) $. Proposition \ref{T:abc} implies that  $\UU_k
 (e)\in\KK$.
 This contradicts   the definition of $e$. Hence
 $\VV(e)=0$.  Thus
 $e=ce_{1,V}$, $|c|=1$ and $ce_{1,V}\in \KK$.

 \begin{remark}
We note that if there is a sequence $n_k\ge1$ such that
$\UU_{n_k}(z_0)$ converges in $H^2$  and $z_0$ satisfies the
hypotheses of Theorem  \ref{T:stab}, then the proof of the
stabilization result obtained in \cite{BCMR} for
finite-dimensional approximations of Schr\"odinger equation works
giving
$$
\UU_{n_k}(z_0,u)\rightarrow e_{1,V}\,\,\,\,\text{in $H^2.$}
$$ However, the existence of such a sequence is an open question.
\end{remark}

\begin{remark} Modifying slightly
   Condition \ref{C:p}, Theorem \ref{T:stab} can be restated for    the  eigenfunction $e_{i,V}$,  $i\ge1$. Indeed,
   one should replace $\la_{1,V}$ and $e_{1,V}$ by $\la_{i,V}$  and
   $e_{i,V}$ in Condition \ref{C:p} and use the   Lyapunov function
\begin{equation*}
\VV_i(z): =\alpha \|( -\Delta+{V} )P_{i,V}z\| ^2+1-|\lag z,e_{i,V}
\rag |^2,\,\,\,\, z\in S\cap H^1_{0}\cap H^2, \end{equation*}where
    $P_{i,V}$ is  the
orthogonal projection in $L^2$ onto the closure of the vector span
of $\{e_{k,V}\}_{k\neq i}$.

\end{remark}

\begin{lemma}\label{L:ner}
Suppose that $r_j\in\R$ and  $r_k\neq r_j$ for $k\neq j$. If
\begin{equation}\label{E:mlml}
\sum_{j=1}^\ty c_je^{ir_jt}=0
\end{equation}
for any  $t\ge0$ and for some sequence  $c_j\in\C$ such that
$\sum_{j=1}^\ty |c_j|<\ty$,   then $c_j=0$ for all $j\ge1$.
\end{lemma}
\begin{proof}Multiplying (\ref{E:mlml}) by $e^{-ir_nt}$ and
  integrating on the interval $[0,T]$, we get
$$
c_n=-\frac{1}{T}\sum_{j=1,j\neq n}^\ty
c_j\int_0^Te^{i(r_j-r_n)t}\dd t= -\frac{1}{T}\sum_{j=1,j\neq
n}^\ty c_j\frac{e^{i(r_j-r_n)T}-1}{i(r_j-r_n)}\rightarrow0
$$as $T\rightarrow\ty$, by the Lebesgue theorem on dominated convergence.
\end{proof}

\subsection{Genericity of Condition \ref{C:p} }\label{S:gen}

Let us recall some definitions. Let $X$ be a complete metric space
and $A\subset X$. Then $A$ is said to be a $G_\delta$ set if it is
a countable intersection of dense open sets. It follows from the
Baire theorem that any $G_\delta$ subset  is dense. A set
$B\subset X$ is called residual if it contains a $G_\delta$
subset.
\begin{example} Let us endow the space  $C^\ty(\overline{ D},\R)$
with  its usual topology given by the countable family of norms:
$$
p_n(Q):=\sum_{|\alpha|\le n}\sup_{x\in D}|\partial^\alpha Q (x)|.
$$
The set ${\cal P} $ of all functions $Q\in C^\ty( \overline{D}
,\R)$ such that $\lag Q e_{1,V},e_{j,V}\rag\neq0$ for all $j\ge1$
is $G_\delta$. Indeed, let us fix an integer $j\ge1$ and let
${\cal P}_j $ be the set of functions $Q\in C^\ty( \overline{D}
,\R)$ verifying $\lag Q e_{1,V},e_{j,V}\rag\neq0$. The unique
continuation theorem for the operator $-\Delta+V$ (see \cite{JK})
implies that there is a ball $B\subset D$ such that
$e_{1,V}(x)e_{j,V}(x)\neq0$ for all $x\in B$. Let $Q\in C^\ty(
\overline{D} ,\R)$ be such that $Q\neq0$, $\supp Q\subset B$ and
$Q\ge0$. Then $Q\in {\cal P}_j$, hence ${\cal P}_j$ is non-empty.
Clearly, ${\cal P}_j$ is open.
 Take any
$Q_1\in  C^\ty( \overline{D} ,\R)$ such that $\lag Q_1
e_{1,V},e_{j,V}\rag=0$ and $Q_2\in{\cal P}_j$. Then $\lag
(Q_1+\tau Q_2)e_{1,V},e_{j,V}\rag\neq0$ for all $\tau\neq0$. Thus
${\cal P}_j$ is dense in $C^\ty( \overline{D} ,\R)$ and ${\cal
P}=\cap_{j=1}^\ty {\cal P}_j$ is a $G_\delta$ set.
\end{example}
The following lemma shows that property (ii) of Condition
\ref{C:p} is generic in 1D case.
\begin{lemma}\label{L:hhh}Let $I\subset\R$ be an interval and let  $\cal Q$ be the set  of all functions $V\in C^\ty(I,\R)$,  verifying
\begin{equation}\label{E:ppo}
\la_{i,V}-\la_{j,V}\neq \la_{p,V}-\la_{q,V}
\end{equation} for all
$i,j,p,q\ge1$ such that
  $\{i,j\}\neq\{p,q\}$ and $i\neq j$. Then $\cal Q$ is a $G_\delta$ set.
\end{lemma}

\begin{proof}
It is well known that the spectrum $\{\la_{j,V}\}$ of
$-\frac{\dd^2}{\dd x^2}+V$ is non-degenerate for any $V\in C^\ty(
\overline{D} ,\R)$, and ${e_{j,V}}$ and ${\la_{j,V}}$ are
real-analytic in $V$ (e.g., see \cite{PT}). Let us introduce the
set  ${\cal Q}_n$, $n\ge1$ of all functions $V\in C^\ty
(\overline{ D},\R)$ such that (\ref{E:ppo}) is satisfied for any
$1\le i,j,p,q\le n$. Clearly,
\begin{align}
{\cal Q}=\bigcap_{n=1}^\ty {\cal Q}_n.\nonumber
\end{align}
It suffices to prove that ${\cal Q}_n$ is open and dense in $C^\ty
(\overline{ D},\R)$. The fact that ${\cal Q}_n$ is open follows
directly from the continuity of ${\la_{j,  V} }$ in $V.$ Let us
prove that  ${\cal Q}_n$ is dense in $C ^\ty(\overline{ D},\R)$.

Take any  $1\le i,j,p,q\le n$ such that
  $\{i,j\}\neq\{p,q\}$ and $i\neq j$, and let ${\cal Q}_{i,j,p,q}$ be the set of
  functions $V\in  C ^\ty(\overline{ D},\R)$ such that (\ref{E:ppo})
  is satisfied.  Suppose we have proved that for any
  $V\in  C ^\ty (\overline{ D},\R)$ there is  $\sigma\in C^\ty
  (\overline{ D},\R)$ such that
\begin{equation}\label{E:kkl}
\la_{i,V+\tau\sigma}-\la_{j,V+\tau\sigma}\neq
\la_{p,V+\tau\sigma}-\la_{q,V+\tau\sigma},
\end{equation}for any small $\tau>0$. This implies that   ${\cal Q}_{i,j,p,q}$
is dense. On the other hand,    ${\cal Q}_{i,j,p,q}$ is open.
Hence ${\cal Q}_n$ is dense, as
$$
 {\cal Q}_n=\bigcap_{1\le i,j,p,q\le n}{\cal
Q}_{i,j,p,q}.
$$

To prove (\ref{E:kkl}), following \cite{ALBERT}, let us write
\begin{align}
\la_{j,V+\tau\sigma}&=\la_{j,V}+\alpha_j\tau+\beta_j(\tau)\tau^2,\label{E:op}\\
e_{j,V+\tau\sigma}&=e_{j,V}+v_j\tau+w_j(\tau)\tau^2.\label{E:op1}
\end{align}
Differentiating the identity
$$(-\frac{\dd^2}{\dd x^2}+V+\tau\sigma-\la_{j,V+\tau\sigma})e_{j,V+\tau\sigma}=0$$
with respect to $\tau$ at $\tau=0$ and using (\ref{E:op}) and
(\ref{E:op1}), we get
$$
(-\frac{\dd^2}{\dd
x^2}+V-\la_{j,V})v_j+(\sigma-\alpha_j)e_{j,V}=0.
$$Taking the scalar product of this identity with $e_{j,V}$, we
obtain
\begin{equation}\label{E:mlk}
 \langle \sigma, (e_{j,V})^2  \rangle=\alpha_j.
\end{equation}Suppose that
$$\la_{i,V+\tau\sigma}-\la_{j,V+\tau\sigma}=\la_{p, V+\tau\sigma}-\la_{q, V+\tau\sigma}$$
for any $\sigma\in C  ^\ty(\overline{ D},\R)$ and for some
sequence $\tau_n\rightarrow0$. Clearly, this implies that
$$
\alpha_i-\alpha_j=\alpha_p-\alpha_q.
$$In view of (\ref{E:mlk}), this gives
\begin{equation}\label{E:aaa}
 (e_{i,V})^2- (e_{j,V})^2= (e_{p,V})^2- (e_{q,V})^2.
\end{equation}On the other hand, by Theorem 9 in \cite{PT} (see page 46), the
system $\{(e_{n,V})^2\}$ is independent for any $V\in L^2$. This
contradiction proves (\ref{E:kkl}) and completes the proof of the
lemma.
\end{proof} We now turn to the multidimensional case.
Let us assume that $D=[0,1]^n$ and introduce the space
\begin{align*} \GG:=\{V\in
C^\ty(D,\R):V(x_1,&\dots,x_n)=V_1(x_1)+\ldots+V_n(x_n)\nonumber\\&\text{for
some $V_k\in C^\ty([0,1],\R),k=1,\ldots,n$}\}. \end{align*} Endow
$\GG$ with the metric of $C^\ty(D,\R)$. It is not difficult to
verify that $\GG$ is a closed subspace in $C^\ty(D,\R)$.
\begin{lemma}The set  of all functions $V\in \GG$,  verifying
\begin{equation}\label{E:ppo2}
\la_{i,V}-\la_{j,V}\neq \la_{p,V}-\la_{q,V}
\end{equation} for all
$i,j,p,q\ge1$ such that
  $\{i,j\}\neq\{p,q\}$ and $i\neq j$,  is a $G_\delta$ set.
\end{lemma}
\begin{proof}
Notice that any eigenfunction of $-\Delta+V$, $V\in\GG$ has the
form
\begin{equation}\label{E:bbb}
e_{l,V}(x_1,\ldots,x_n)=e_{l_1,V_1}(x_1)\cdot\ldots\cdot
e_{l_n,V_n}(x_n),\end{equation}where $e_{l_k,V_k}(x_k)$ is an
eigenfunction of the operator $-\frac{\dd^2}{\dd x_{k}^2}+V_k$.
Indeed, any function of   form (\ref{E:bbb}) is an eigenfunction,
and the set of all functions of this form is a basis in $L^2(D)$.

Let $i,j,p,q\ge1$ be such that
  $\{i,j\}\neq\{p,q\}$ and $i\neq j$,  and let $e_{i_n,V_n}(x_n)  $, $ e_{j_n,V_n}(x_n)  $,
$ e_{p_n,V_n}(x_n)  $ and $ e_{q_n,V_n}(x_n)  $ be the
eigenfunctions in (\ref{E:bbb}). Without loss of generality, we
can suppose that the functions $(e_{i_n,V_n}(x_n))^2$,
$(e_{j_n,V_n}(x_n))^2$, $(e_{p_n,V_n}(x_n))^2$ and
$(e_{q_n,V_n}(x_n))^2$
   are linearly independent (see Theorem
9 in \cite{PT}).
 Any eigenfunction $e_{l_k,V_k}$ has a finite number of zeros in
interval  $[0,1]$. Hence, choosing appropriately the point $x^*\in
[0,1]^{n-1}$, we see that the functions $(e_{i,V}(x^*,x_n))^2$,
$(e_{j,V}(x^*,x_n))^2$, $(e_{p,V}(x^*,x_n))^2$ and
$(e_{q,V}(x^*,x_n))^2$, $x_n\in[0,1]$ are linearly independent.
This implies that relation (\ref{E:aaa}) does not hold. Thus the
proof of Lemma \ref{L:hhh} works  implying the genericity.
\end{proof}

\section{Applications}\label{S:Application}

\subsection{Nonlinear  Schr\"odinger equation }\label{S:Application} Let us  consider
the nonlinear Schr\"odinger equation
\begin{align}
i\dot z   &= -\Delta z+V(x)z+ u(t) Q(x)|z|^2z,\,\,\,\,x\in D,\label{E:hav21}\\
z\arrowvert_{\partial D}&=0,\label{E:ep21}\\
 z(0,x)&=z_0(x),\label{E:sp21}
\end{align}where $D\subset\R^3$ is a bounded domain with smooth
boundary. Problem (\ref{E:hav21})-(\ref{E:sp21}) is locally
well-posed.
\begin{lemma}\label{L:LD2}

For any   $z_0\in H_0^1\cap H^2$   and for any
  $u\in L^1_{loc}([0,\ty),\R)$ there is a time $T>0$ such that
  problem
(\ref{E:hav21})-(\ref{E:sp21}) has a unique solution $z\in
C([0,T],H^2)$. Furthermore, the resolving operator
$\UU_t(\cdot,u):H_0^1\cap H^2 \rightarrow H_0^1\cap H^2 $ taking
$z_0$ to
    $z(t)$ satisfies the relation
\begin{align*}
\|\UU_t(z_0,u)\|&=\|z_0\|,\,\,\,t\in[0,T].
\end{align*}
\end{lemma} See \cite{CW} for the proof. Define
$z(t)=\UU_t(z_0,u)$ and   let us calculate the
derivative\begin{align*} \frac{\dd}{\dd
t}\VV(z(t))&=2\alpha\Re(\langle (-\Delta +V) P_{i,V}{\dot{
z}},P_{i,V} (-\Delta +V)  z\rangle)-2\Re(\lag\dot z,e_{i,V}\rag\lag e_{i,V},z\rag)\nonumber\\
&= 2\alpha\Re(\langle  (-\Delta +V) P_{i,V}(i\Delta z -iVz-iu
Q|z|^2z ), (-\Delta
+V)  P_{i,V}z\rangle)\nonumber\\
&\quad-2\Re(\lag  i\Delta z -iVz-iuQ|z|^2z,e_{i,V}\rag\lag e_{i,V},z\rag)\nonumber\\
&= 2u\Im(\alpha\lag (-\Delta +V)  P_{i,V}(Q|z|^2z),
(-\Delta +V) P_{i,V}z\rangle\nonumber\\
&\quad-\lag Q|z|^2z,e_{i,V}\rag\lag e_{i,V},z\rag).
\end{align*}Take
\begin{equation}\label{E:vhh}
  {{ {u}}}(z):=-\Im(\langle \alpha (-\Delta +V)   P_{i,V}(Q|z|^2z ), (-\Delta
 +V)
 P_{i,V}z\rangle
 - \lag Q|z|^2z,e_{i,V}\rag\lag e_{i,V},z\rag).
\end{equation}
Problem (\ref{E:hav21})-(\ref{E:sp21}) with feedback (\ref{E:vhh})
is globally well-posed in $H^2$ (cf. Theorem \ref{T:abc}).  Let
$\UU_t(\cdot): H_0^1\cap H^2\rightarrow
   H_0^1\cap H^2$ be the resolving operator. To
formulate the main result, we introduce the following hypothesis.
\begin{condition}\label{C:2p}
The functions $V,Q\in C^\ty(\overline{D},\R) $ are such that:
\begin{enumerate}
\item [(i)] $\langle Qe_{i,V}e_{j,V},e_{p,V}e_{q,V}\rangle\neq0$
  for all $i, j,p,q\ge1$,
\item  [(ii)]$\la_{i,V}-\la_{j,V}+\la_{p,V}-\la_{q,V}\neq \la_{i',V}-\la_{j',V}+\la_{p',V}-\la_{q',V}$ for all  integers $ i,j,p,q,i',j',p',q'$ such that
  $\{i,j,p,q\}\neq\{i',j',p',q'\}$ and $\{i,p\}\neq\{j,q\}$.

\end{enumerate}
\end{condition}

The  theorem below is the version of Theorem \ref{T:stab} for
system (\ref{E:hav21})-(\ref{E:sp21}).

\begin{theorem}\label{T:2mm}
 Under Condition \ref{C:2p},  there is a finite  or countable set $J\subset \R^*_+$ such that for any
  $\alpha\notin J$, $\l\ge1$ and $z_0\in S\cap H^1_{0}\cap H^2$ with $\lag z_0,e_{l,V} \rag\neq0$
  and $0<~\VV_l(z_0)<1$  there is a sequence
 $k_n\ge1$ verifying
 $$
 \UU_{k_n}(z_0 )\rightharpoonup  c e_{l,V}\,\,\,\,\text{in
$H^2$},
 $$where $c\in\C,|c|=1.$ \end{theorem}

The proof of this theorem is very close to that of Theorem
\ref{T:stab}. One should notice that, under Condition \ref{C:2p},
there is a time $t_0>0$ such that $ {{
{u}}}(\UU_{t_0}(z_0,0))\neq~0$   and then conclude as in Step 2 in
the proof of   Theorem \ref{T:stab}.
\begin{remark}Notice that, as equation (\ref{E:hav21}) is
nonlinear, the distance between two solutions with the same
control is not constant. Hence the proof of approximate
controllability given in  Theorem \ref{T:con} does not work here.

\end{remark} \begin{lemma}
 For any $l\ge1$, $d>0$ and  $z_0\in S$ there is a control $u\in C_d$
 and a time $k\ge1$ such that  $\lag\UU_k(z_0,u),e_{l,V}\rag\neq0.$
 \end{lemma}
 \begin{proof}  Suppose that $\lag z_0,e_{l,V} \rag=0$. Let us show that there is a  control $u\in
 C_d$ such that    $\lag \UU_k (z_0,u),e_{l,V} \rag\neq0$ for some
 $k\ge0$. As (\ref{E:hav21})  is nonlinear, the proof given in
 Lemma \ref{L:eee} does not work.

If $z_0\notin\{ ce_{j,V}: c\in\C,|c|=1, j\ge1\}$, then, by Theorem
\ref{T:2mm},  there is an integer $p\ge1$,    sequence  $k_n\ge1$
and   constant $c\in~\C, |c|=1$ such that $\UU_{k_n} (z_0 )
\rightharpoonup ce_{p,V} $ in $H^{2}$. Hence, without loss of
generality, we can suppose that $z_0=e_{p,V}$ for some $p\neq l$.
Let us introduce the following two-dimensional subspace of
$L^2([0,1],\R)$:
$$
E=\{a\sin (\la_{p,V}-\la_{l,V}) t+b\cos(\la_{p,V}-\la_{l,V}) t:
a,b\in\R   \}.
$$
For any   $u\in E$, define the mapping
$\Phi(u)=\lag\UU_1(e_{p,V},u), e_{l,V}\rag$, whenever the solution
$\UU_t(e_{p,V}, u)$ exists up to time $t=1$. Notice that $\Phi(
0)=\lag e^{-i\la_{p,V}}e_{p,V},e_{l,V}\rag=0$, hence $\Phi $ is
well defined in a neighborhood of $0\in E$. We are going to show
that the conditions of inverse mapping theorem are satisfied in a
neighborhood of the point $0\in   E $. Clearly, $\Phi$ is
continuously differentiable. Let us show that mapping $D\Phi ( 0):
E \rightarrow \C$ is an isomorphism.
 Consider the
linearization of (\ref{E:hav21}), (\ref{E:ep21}), $z_0=e_{p,V}$
around $(e^{-i \la_{p,V} t}e_{p,V},0)$:
\begin{align}
i\dot y   &= -\Delta y+V(x)y+ u(t)
Q(x)e_{p,V}^3e^{-i\la_{p,V}t},\,\,\,\,x\in D.\label{E:lin}\\
y\arrowvert_{\partial D}&=0,\label{E:ep41}\\
 y(0)&=0.\label{E:sp41}
\end{align} One can verify that
$D\Phi( 0)( u)=\lag y(1),e_{l,V}\rag$.  System
(\ref{E:lin})-(\ref{E:sp41})    is equivalent to
\begin{equation}\label{E:D}
y=  -i\int_0^t e^{-i\la_{p,V}s}u(s)S(t-s)(Q e_{p,V}^3)\dd s,
\end{equation}where $S(t)$ is the unitary group associated with
  $i\Delta-iV.$ Taking the scalar product of (\ref{E:D}) with
$ e_{l,V}$, we obtain for $t=1$
 \begin{equation*}
\lag y,e_{l,V}\rag= -ie^{-i\la_{l,V}} \lag  Q e_{p,V}^3 ,
e_{l,V}\rag\int_0^1 e^{-i(\la_{p,V}-\la_{l,V})s}u(s)\dd   s.
\end{equation*} Condition \ref{C:2p} implies that $\la_{p,V}-\la_{l,V}\neq
0$, hence $D\Phi( 0):E\rightarrow\C$ is an isomorphism. Applying
the inverse mapping theorem, we conclude that $\Phi$ is  $C^1$
diffeomorphism in a neighborhood of $ 0\in E$. Thus there is a
control $u\in C_d$ such that    $\lag \UU_1 (e_{p,V},u),e_{l,V}
\rag\neq0.$
\end{proof}

\subsection{Randomly forced Schr\"odinger equation }\label{intro}
\subsubsection{Growth of Sobolev norms}\label{intr} Let us consider the
problem
\begin{align}
i\dot z &= -  \Delta z+V(x)z +\beta(t) Q(x)z,\,\,\,\,x\in D,\label{E:havvv}\\
z\arrowvert_{\partial D}&=0,\label{E:epvv}\\
 z(0)&=z_0,\label{E:spvv}
\end{align}where   $V,Q\in C^\ty(\overline{D},\R)$ are   given
functions. We  assume that $\beta(t)$ is a random process of the
form (\ref{E:pu}), where
  the random variables $\eta_k$  verify the following condition.

\begin{condition}\label{C:c}
The random variables $\eta_k$  have the form
$$
\eta_k(t)=\sum_{j=1}^\infty b_j\xi_{jk}g_j(t), \,\,\,\, t\in[0,1],
$$where $\{g_j\}$ is an orthonormal basis in $L^2([0,1],\R)$,
$b_j>0$ are constants with
$$
 \sum_{j=1}^\infty b_j^2<\infty,
$$ and $\xi_{jk}$ are independent real-valued random variables such that
$\E\xi_{jk}^2=1$. Moreover, the distribution of $\xi_{jk}$
possesses a continuous density $\rho_j$ with respect to the
Lebesgue measure and $\rho_j(r)>0$ for all $r\in\R$.
\end{condition}

Notice that  this condition in particular implies that
$$
\pP\{\|u-\beta\|_{L^2([0,l])}<\e\}>0
$$for any $u\in L^2([0,l])$
and $\e>0$. Moreover, using the continuity of the mapping
$\UU_l(z_0,\cdot):L^2([0,l])\rightarrow L^2(D)$, for any $\de>0$
we can find   a constant $\e>0$ such that
$$
\pP\{\|\UU_l(z_0,\beta)-\UU_l(z_0,u)\|<\delta\}\ge\pP\{\|u-\beta\|_{L^2([0,l])}<\e\}>0.
$$ Hence, any
point $\UU_l(z_0,u), u\in L^2([0,l])$ is in the support of the
measure $\DD(\UU_l(z_0,\beta))$.

 The following theorem is the main result of this section.
\begin{theorem}\label{T:him}
Suppose that Conditions \ref{C:p} and \ref{C:c} are satisfied.
Then   for any $s>0$ and $z\in H^s \backslash\{0\}$ we have
\begin{equation}\label{E:hhh1}
\pP\{\limsup_{k\rightarrow\infty}\|\UU_k(z,\beta)\|_s=\infty\}=1.
\end{equation}
\end{theorem}
\subsubsection{Proof of Theorem \ref{T:him}}\label{S:him} By Theorem \ref{T:con},   system (\ref{E:hav1}),
 (\ref{E:ep1}) is approximately controllable at integer times.
 Since the equation is linear in $z$, 
  it suffices to prove (\ref{E:hhh1}) for  any $z\in
 S\cap H^s $. Without loss of generality, we can assume that  $s\in (0,2]$.

\vspace{6pt} \textbf{Step 1.} Let us fix a constant $r>0$ and
introduce the stopping time
$$
\tau_r(z)=\min\{k\ge0:\UU_k(z,\beta)\in
B_{H^{-s}}(0,r)\},\,\,\,\,z\in
 B_{L^2}(0,1).
$$
Then we have
\begin{equation}\label{E:hay}
\pP \{\tau_r(z)<\infty\}=1.
\end{equation} Indeed, choose an
arbitrary point $z'\in  S\cap B_{H^{-s}}(0,r)$. By the property of
approximate controllability in $L^2$, there is a control $u\in
C_d$ such that $\UU_l(z,u)$ is sufficiently close to $z'$ in
$L^2$, hence $\UU_l(z,u)\in B_{H^{-s}}(0,r)$. As $\UU_l(z,u)$ is
in the support of measure $\DD(\UU_l(z ,\beta))$, we have
$$
\pP\{\UU_l(z,\beta)\in B_{H^{-s}}(0,r)\}>0.
$$Using the continuity of the resolving operator in negative Sobolev
norms, we see that there is  an $H^{-s}$-neighborhood $\OO=\OO(z)$
of $z$ such that
\begin{equation}
 \sup_{{y\in \OO}}\pP \{\tau_r(y)
>l\}<1.\nonumber
\end{equation}
From the compactness of $B_{L^2}(0,1)$ in $H^{-s}$ it follows that
there is a  time $k\ge1$ such that
\begin{equation}\label{E:vah}
a:=\sup_{y\in B_{L^2}(0,1)} \pP \{\tau_r(y)
>k\}<1.
\end{equation}
Using the Markov property  and (\ref{E:vah}), we obtain
\begin{align}
\pP\{\tau_r(y) >nk\}&=\E (I_{\{\tau_r(y) >(n-1 )k\}}\pP\{\tau_r(x)
>k\}|_{x= {\UU_{(n-1)k}(y,\beta)}})\nonumber\\&\le a\pP \{\tau_r(y)
>(n-1)k\}.
\end{align} Hence
\begin{equation}
\pP \{\tau_r(y) >nk\}\le a^n.\nonumber
\end{equation}
Using the Borel--Cantelli lemma, we arrive at  (\ref{E:hay}).

\vspace{6pt}  \textbf{Step 2.} Take any $z\in S\cap H^s $.
 Choosing   $r=\frac{1}{n}$  and using (\ref{E:hay}), we get
 \begin{equation}\label{E:hhh2}
\pP\{\liminf_{k\rightarrow\infty}\|\UU_k(z,\beta)\|_{-s}=0\}=1.
\end{equation} Define the event
$$
\aA:=\{\omega\in\Omega:\limsup_{k\rightarrow\infty}\|\UU_k(z,\beta)\|_s<\infty\}.
$$
Suppose  that
\begin{equation}
\pP\{\aA\}>0.\nonumber
\end{equation}
By (\ref{E:hhh2}), for almost any $\omega\in\aA$ there is a
sequence $n_k\rightarrow\ty$ such that
\begin{equation}\label{E:Pp}
\lim_{n\rightarrow\ty}\|\UU_{n_k}(z,\beta)\|_{-s}=0.
\end{equation}On the other hand, for any $\omega\in\aA$, there is a subsequence of
$n_k$ (which is also denoted by $n_k$) and an element $w\in S$
such that
$$
\|\UU_{n_k}(z,\beta)-w\|\rightarrow0.
$$
This contradicts (\ref{E:Pp}). Thus $\pP\{\aA\}=0$.

\begin{remark}In view of Theorem \ref{T:2mm}, under Condition \ref{C:2p}, Theorem \ref{T:him}
holds also in the  case of nonlinear equation (\ref{E:hav21}). The
proof is literally the same. One should just pay attention to the
fact that, as in this case
       finite time blow-up is possible, the restriction of the solution
     at integer times forms a Markov chain with values in $H^s\cup
     \{\ty\}$ (e.g., see \cite{REV}).

\end{remark}

\end{document}